\documentclass[12pt,oneside]{article}
\usepackage{amssymb,amsfonts,amsmath,epsfig}
\usepackage{natbib}
\usepackage{hyperref}

\numberwithin{equation}{section}
\newtheorem{theorem}{Theorem}[section]
\newtheorem{prop}{Proposition}[section]
\newtheorem{lemma}{Lemma}[section]
\newtheorem{cor}{Corollary}[section]

\def\eqref#1{(\ref{#1})}
\def\vf#1{{\boldsymbol{#1}}}
\def\var{\mathrm{Var}}

\def\Grp#1{\left(#1\right)}
\def\Cbr#1{\left\{#1\right\}}
\def\Sbr#1{\left[#1\right]}

\def\Abs#1{\left|#1\right|}

\def\Norm#1{\left\|#1\right\|}
\def\cf#1{\mathbf{1}\Cbr{#1}}
\def\Pr{\mathbb{P}}
\def\pfdr{{\rm pFDR}}
\def\power{{\rm power}}

\def\Limsup{\mathop{\overline{\mathrm{lim}}}}
\def\eno#1#2{#1_1, \ldots, #1_#2}

\def\inum#1{{#1}_1, {#1}_2, \ldots}

\def\Reals{\mathbb{R}}
\def\qed{\hfill$\Box$\medbreak}
\def\gv{\,|\,}

\def\nth#1{\frac{1}{#1}}

\def\Sp#1{\sp{(#1)}}

\def\pa{a} 

\def\convP{\stackrel{\rm P}{\to}}
\def\toi{\to\infty}

\def\dgamma{\text{Gamma}}
\def\cE{\mathcal{E}}
\def\shape{\nu}

\long\def\comment#1{}

\begin{document}

\noindent
DATA VOLUME AND POWER OF MULTIPLE TESTS WITH SMALL SAMPLE
SIZE PER NULL

\vskip 3mm

\vskip 8mm
\noindent Zhiyi Chi\footnote{Research partially supported by NSF
  DMS~0706048 and NIH MH~68028.}

\noindent Department of Statistics

\noindent University of Connecticut

\noindent 215 Glenbrook Road, U-4120

\noindent Storrs, CT 06269, USA

\noindent zchi@stat.uconn.edu


\vskip 3mm
\noindent ABSTRACT

In multiple hypothesis testing, the volume of data, defined as the
number of replications per null times the total number of nulls,
usually defines the amount of resource required.  On the other hand,
power is an important measure of performance for multiple testing.
Due to practical constraints, the number of replications per null may
not be large enough in terms of the difference between false and true
nulls.  For the case where the population fraction of false nulls is
constant, we show that, as the difference between false and true nulls
becomes increasingly subtle while the number of replications per null
cannot increase fast enough, (1) in order to have enough chance that
the data to be collected will yield some trustworthy findings, as
measured by a conditional version of the positive false discovery rate
(pFDR), the volume of data has to grow at a rate much faster than in
the case where the number of replications per null can be large
enough, and (2) in order to control the pFDR asymptotically, power has
to decay to 0 in a rate highly sensitive to rejection criterion and
there is no asymptotically most powerful procedures among those that
control the pFDR asymptotically at the same level.

\section{Introduction} \label{sec:intro}
Multiple hypothesis testing often faces situations where distributions
under false nulls only have a subtle difference from those under true
nulls.  To make the difference evident, it is necessary to make
repeated measurements from the distributions.  As is well known, in
order to attain a fixed power, the number of replications for each
null, henceforth denoted by $k$, should be large enough.  Roughly
speaking, if the difference between the distributions under false
nulls and those under true nulls is $\delta$, then $k$ should be of
the same order as $\delta^{-2}$.  However, in practice, oftentimes $k$
cannot be large enough.  There are many reasons for this: the time
window that allows measurements is finite, the experimental unit
associated with each null can only sustain a limited number of
measurements, and so on.  Under this circumstance, there are two
important and related issues.

First, if the underlying objective of the tests is to identify even
just a few false nulls irrespective of power, so that one can get
useful clues for further study, then in some cases it may be reached
by testing a large number of nulls.  For example, suppose a population
has a small fraction of ``atypical'' units.  Then, in order to capture
at least one of them, an approach is to obtain a large sample from the
population, and, for each sampled unit, determine whether to reject
the null that it is typical.  In this case, even though at population
level, the difference between atypical and typical units may be small,
there is a chance that some atypical units will ``show up'' with
pronounced differences from typical ones, making them easy to detect.
In order to increase the chance, one would hope $N$, the number of
examined units, to be as large as possible.  Oftentimes, however, the
amount of resources required for the tests is in proportion to $N$ or
$V=kN$, the ``volume'' of data.  This imposes a constraint on $N$ and
$V$ and raises the following question: provided that $k$ cannot be
large enough, what is the minimum value of $N$ or $V$ in order to have
desirable test results?

The answer to the question depends on what performance criterion to
use for the tests.  A useful criterion is $\pfdr \le \alpha$, where
$\alpha\in (0,1)$ is a pre-specified level and $\pfdr=E[R_0/R\gv
R>0]$ is known as the positive false discovery rate, with $R$ the
number of rejected nulls and $R_0$ that of rejected true nulls
\citep{storey:03}.   Comparing to the false discovery rate
$E[R_0/(R\vee 1)]$ \citep{benjamini:hoc:95}, the pFDR is more useful
when the objective is to reject \emph{some\/} nulls \citep{chi:07,
  chi:tan:08}.  In view of this, it is desirable to apply the idea
of pFDR to data directly.  Denote by $\vf X$ the data \emph{to be\/}
collected for the tests.  Under a Bayesian 
framework, we propose the following variant of the pFDR criterion,
\begin{align}
  \Pr\Grp{
    \inf E\Sbr{R_0/R\gv R>0,\, \vf X}\le \alpha
  }\ge p,
  \label{eq:criterion}
\end{align}
where $E[R_0/R\gv R>0, \vf X]$ is defined to be 1 if $R=0$, $\alpha$,
$p\in (0,1)$ are pre-specified constants, and the infimum is taken
over all applicable multiple testing procedures that are solely based
on the data; cf.\ \eqref{eq:multitest}.  For any such procedure, once
$\vf X$ is given, $R$ is determined, whereas $R_0$ remains random,
with its distribution determined by the posterior likelihoods of the
truths of the nulls.  This gives rise to the conditional expectation.
As a result, the infimum in \eqref{eq:criterion} is a function of $\vf
X$.  The criterion means that, with probability no less than $p$, 
the data to be collected will yield one or more rejections, which,
under a future study for verification, have an expected
fraction of false rejections no more than $\alpha$.  The minimum value
of $N$ or $V$ will be studied under this criterion.

The second issue is straightforward, that is, even when an arbitrarily
large number of nulls can be tested, power may still be an important
concern.  Provided $k$ cannot be large enough, how much power can be
hoped for?

The issue of minimum $N$ or $V$ can be thought of as a dual to the
issue that given $N$, how small the differences between false 
and true nulls can be before it becomes virtually impossible to detect
false nulls; see \cite{donoho:jin:04} and references therein.  The
issue of power is more extensively studied; see \cite{efron:07} and
references therein.  However, the setup here is different.  First,
both issues are considered in relation to $k$ and it is necessary to
take into account the fact that the distributions of test statistics
not only depend on the underlying data generating distributions, but
also on $k$, the number of replications from the distributions for
each null.  Second, no sparsity is assumed for false nulls.  Instead,
false nulls are assumed to be increasingly similar to true nulls,
while $k$ cannot increase fast enough to compensate for this.  As a
result, the test statistics provide increasingly weak evidence to
separate false nulls from true nulls.

All our results are obtained for the case where distributions under
false nulls and those under true nulls are known and the population
fraction of false nulls is known as well.  In practice, especially
in exploratory studies, while there may be relatively good knowledge
about distributions under true nulls, oftentimes there is little
knowledge about distributions under false nulls or the population
fraction of false nulls.  Thus, the case we consider is an ideal
one and the results provide limits on what can be achieved in more
realistic cases.

Section \ref{sec:setup} covers preliminaries and identifies the
quantity central to the analysis.  Main results are stated in Section
\ref{sec:main}.  Section
\ref{subsec:lrt} considers the asymptotics of the minimum $N$ and $V$
for general parametric models.  Using Cram\'er-type large deviations,
the 
results are established for the case where the growth of $k$ is just a
little slower than $\delta^{-2}$.  It shows that in this case, the
minimum $N$ and $V$ have to grow much faster  than in the case where
$k$ is of the same order as $\delta^{-2}$.  Sections
\ref{subsec:normal} and \ref{subsec:gamma} obtain refined results for
tests on mean values of normal distributions with known variance and
scaled gamma distributions.  Section \ref{subsec:power} considers the
power of multiple tests as $k$ can not increase as rapidly as
$\delta^{-2}$ as $\delta\to 0$.  It shows that for procedures that
asymptotically control the pFDR at a given level,
the power decreases to 0 and is highly sensitive to small changes in
rejection criterion, and consequently, there is no asymptotically 
most powerful procedure among such procedures.  Section
\ref{sec:summary} concludes the article with a summary and remarks.
The proofs of the main results are collected in the Appendix.

\section{Preliminaries} \label{sec:setup}
Denote by $\eno H N$ the nulls.  Suppose that for each $H_i$, a
sample $X_{i1}$, \ldots, $X_{i k}$ is collected.  Let $\eta_i = \cf{H_i
  \text{ is false}}$.  The analysis will be under the
following random effects model:
\begin{gather}
  \begin{array}{c}
    (\eta_i, X_{i1}, \ldots, X_{i k}), \ \ i\ge 1, \text{ are
      i.i.d.\ such that}\\
    \Pr(\eta_i = 1) = \pa \text{ and } \\
    \begin{array}{ll}
      \text{given $\eta_i=0$},&
      X_{i1}, \ldots, X_{i k} \ \text{ are i.i.d.}\ \sim f_0(x), \\
      \text{given $\eta_i=1$},&
      X_{i1}, \ldots, X_{i k} \ \text{ are i.i.d.}\ \sim f_a(x),
    \end{array}
  \end{array}
  \label{eq:model}
\end{gather}
where $0<\pa<1$ is the population fraction of false nulls, and $f_0$,
$f_a$ are probability densities under true nulls and false nulls,
respectively \citep{efron:etal:01, genovese:was:02}.  In this article,
we will only consider the case where $\pa$ is fixed.

Following \cite{chi:tan:08}, a multiple testing procedure is a
deterministic mapping
\begin{align}
  d(\vf X) = (d_1(\vf X), \ldots, d_N(\vf X)),  \label{eq:multitest}
\end{align}
such that $H_i$ is rejected if and only if $d_i(\vf X)=1$.  Under the
random effects model, $\vf X=\{X_{kl}\}$ and it can be shown
that the criterion \eqref{eq:criterion} can be rewritten as
\begin{gather}
  \Pr(\text{at least one $\cE_i(k)$ occurs},\ i=1,\ldots, N)\ge p,
  \label{eq:detect} \\
  \text{where}\quad
  \cE_i(k):=
  \Cbr{X_{i1}, \ldots, X_{ik} \text{ satisfy }
    \prod_{j=1}^k \frac{f_a(X_{i j})}{f_0(X_{i j})} \ge 
    \frac{(1-\pa) (1-\alpha)}{\pa\alpha}
  };
  \label{eq:pfdr}
\end{gather}
see Section \ref{sec:proofs} for a sketch of proof.

For fixed $k$, the events $\cE_i(k)$ are independent of each other and
have the same probability, which depends on both $k$ and
the difference between $f_0$ and $f_a$.  Suppose the difference
can be parameterized by $\delta$.  For example, if $f_0=N(\theta_0,1)$ and
$f_a=N(\theta, 1)$, then $\delta$ can be taken as $\theta-\theta_0$.
Denote the common probability of $\cE_i(k)$ by
\begin{align}
  p_{k,\delta}(\alpha) = \Pr(\cE_i(k)), \quad i=1,\ldots, N.
  \label{eq:p-def}
\end{align}
Then \eqref{eq:detect} is equivalent to $1 -
(1-p_{k,\delta}(\alpha))^N \ge p$, or $N\ge \ln (1-p)/\ln
(1-p_{k,\delta}(\alpha))$.

The case we will focus on is where the false nulls are increasingly
similar to true nulls while $k$ cannot increase fast enough to
compensate for the decreasing $\delta$; more specifically, $\delta\to
0$ while $k\to\infty$ at a slower rate than $\delta^{-2}$.  Then, as
maybe expected, $p_{k,\delta}(\alpha)\to 0$.  Provided that the growth
of $k$ makes sure $p_{k,\delta}(\alpha)>0$, the minimum number of
nulls and volume of data to satisfy \eqref{eq:criterion} are
\begin{align}
  N_* = \frac{(1+o(1))}{p_{k,\delta}(\alpha)} \ln\nth{1-p}, \quad
  V_* = k N_* = \frac{(1+o(1)) k}{p_{k,\delta}(\alpha)} \ln\nth{1-p},
  \label{eq:min-n}
\end{align}
respectively.  Thus, the main task of the analysis is to find the
asymptotic of $p_{k,\delta}(\alpha)$.  Note that under the random
effects model,
\begin{align}
  p_{k,\delta}(\alpha)= (1-a) \Pr_0(\cE_i(k))
  + a \Pr_a(\cE_i(k)),  \label{eq:pfdr2}
\end{align}
where $\Pr_0$ is the probability measure under $f_0$ and $\Pr_a$ that
under $f_a$.

Finally, some comments on the random effects model.  It may be
desirable to relax the assumption that the data distributions under
false nulls are identical.  To do this, one approach is to assume that
under false nulls, the data obeys another random effects model, such
that given $\eta_i=1$, a parameter $\theta_i$ is first drawn from a
distribution $G$, and then 
$X_{i1}$, \ldots, $X_{i k}$ are drawn from $f_{\theta_i}\not=f_0$
\citep{genovese:was:02}.  However, by letting $f_a = \int
f_\theta\,dG(\theta)$, it is seen that this model can be treated in
the same way as \eqref{eq:model}.  Another approach is to use Poisson
approximation to evaluate the probability in \eqref{eq:detect}, which
does not require the distributions under false nulls be identical, and
can even allow weak dependency between $X_{ij}$ \citep{arratia:90a}.
However, a full development of the Poisson approach is beyond the
scope of the paper.

\comment{Even without 
assumptions on random effects, the data \emph{directly\/} employed
by the tests may still follow identical distributions under false
nulls.  For example, consider $t$-tests for $H_i: \mu_i=0$ in $N(\mu_i,
\sigma_i^2)$, where nothing is known about $\sigma_i^2$, and each
$t$-test is based on an i.i.d.\ sample of $k$ observations from $N(\mu_i,
\sigma_i^2)$.  If $\mu_i = c\sigma_i$ for false nulls, where
$c>0$ is an unknown constant, then it is seen that the $t$-statistics
obey \eqref{eq:model} with $f_0=t_{k-1}$ and $f_a=t_{k-1,d}$, where
$t_{k-1,d}$ is the noncentral $t$-distribution with $k$ degrees of
freedom and noncentrality $d=\sqrt{k}c$.}

\section{Statement of main results} \label{sec:main}
Recall that $\pa$ is assumed to be fixed.  Henceforth, denote
$Q_\alpha = (1/\pa-1)(1/\alpha-1)$.

\subsection{Data volume for general multiple tests}
\label{subsec:lrt}
Suppose the observations $X_{i j}$ take values in a Euclidean space
$\Omega$ and both $f_0$ and $f_a$ belong to a parametric family of
densities $\{f_\theta, \theta\in \Theta\}$ with respect to the
Lebesgue measure $dx$ on $\Omega$, where $\Theta$ is an open set in
$\Reals^d$.  Suppose $f_0 = f_{\theta_0}$ and $f_a = f_\theta$.  Let
$\ell(\theta,x) = \ln f_\theta(x)$.  Suppose that for each
$x\in\Omega$, $\ell(\theta, x)$ is twice differentiable with
\begin{align*}
  \dot\ell(\theta,x) =
  \Sbr{\frac{\partial\ell(\theta,x)}{\partial\theta_1},
    \ldots, \frac{\partial\ell(\theta,x)}{\partial\theta_d}}^T, \ \
  \ddot\ell(\theta,x) =
  \Sbr{\frac{\partial^2\ell(\theta,x)}
    {\partial\theta_k\partial\theta_l}}.
\end{align*}
We assume that $\{f_\theta\}$ satisfies regular conditions so that
\begin{gather*}
  E_\theta[\dot\ell(\theta, X)] = 0, \quad
  \var_\theta[\dot\ell(\theta, X)] = -E_\theta[\ddot\ell(\theta,X)] =
  I(\theta),
\end{gather*}
where $I(\theta)$ is the Fisher information and $E_\theta$ and
$\var_\theta$ denote the expectation and variance under $f_\theta$,
respectively.  For each $\theta, \theta'\in \Theta$, by Taylor
expansion,
\begin{align*}
  \ell(\theta', x) - \ell(\theta, x)
  &
  =
  \dot\ell(\theta, x)^T(\theta' - \theta)
  + \frac{(\theta'-\theta)^T \ddot\ell(\theta,x)(\theta'-\theta)}{2}
  + o(|\theta'-\theta|^2) \\
  &
  =
  \dot\ell(\theta, x)^T(\theta' - \theta)
  + \frac{(\theta'-\theta)^T A(\theta, \theta', x)
    (\theta'-\theta)}{2},
\end{align*}
where $A(\theta,\theta',x)$ is a $d\times d$ symmetric matrix.  Under
regular conditions, one would expect that as $\theta'\to\theta$, 
$E_\theta [A(\theta, \theta', X)] \to E_\theta[\ddot\ell(\theta, X)] =
-I(\theta)$.  However, for the analysis, a few stronger assumptions
are needed.

\paragraph{Assumptions.}
\begin{enumerate}
\item \label{A:I-pd} $I(\theta_0)$ is positive definite.
\item \label{A:I-diff} There are $C>0$ and $r>0$, such that
  $|I(\theta)-I(\theta_0)|\le C|\theta - \theta_0|$ if
  $|\theta-\theta_0| < r$.
\item \label{A:mnt} $M_3:=\sup_{\theta\in \Theta} E_\theta
  |\dot\ell(\theta, X)|^3<\infty$.
\item \label{A:ldp} For any $\epsilon>0$, there are positive numbers
  $r$, $\lambda_0$ and $\lambda$, such that for any $\theta$ and
  $\theta'\in \Theta$ with $|\theta-\theta'| < r$,  if $\inum X$ are
  i.i.d.\ $\sim f_\theta$, then
  \begin{align}
    \Pr\Grp{
      \Norm{
        \nth k\sum_{i=1}^k A(\theta, \theta', X_i) + I(\theta) 
      } > \epsilon
    } \le \lambda_0 e^{-\lambda k}
    \label{eq:LDP}
  \end{align}
  for all $k\ge 1$, where for any matrix $M=(m_{i j})$, $\|M\|=\max
  |m_{i j}|$. 
\end{enumerate}

Assumption 1 is standard.  Assumption 2 holds if $I(\theta)$
is differentiable at $\theta_0$, which along with Assumption 3 is
satisfied by many parametric models.  Assumption 4 is not hard to
verify by using the fact that the normed quantity in \eqref{eq:LDP} is
bounded by $D_1+D_2$, where
$$
D_1=\Norm{h(\theta, \theta') + I(\theta)}, 
\quad
D_2=
\Norm{
  \nth k\sum_{i=1}^k A(\theta, \theta', X_i) - h(\theta, \theta')
},
$$
with $h(\theta, \theta')=E_\theta[A(\theta, \theta', X)]$.  In fact,
since $h(\theta, \theta) = -I(\theta)$, provided that $h(\theta,
\theta')$ is uniformly continuous, $D_1$ is uniformly small for 
$|\theta-\theta'|\ll 1$.  On the other hand, exponential inequalities
can be used for $D_2$.  For instance, if $|A(\theta, \theta', X)|\le M$
for some nonrandom $M>0$ for all $\theta$, $\theta'$ and $X\sim
f_\theta$, then Hoeffding's inequality gives $\Pr(D_2\ge \epsilon) \le
\lambda_0 e^{-\lambda k}$ for some $\lambda_0$, $\lambda>0$
\citep{pollard:84}.  As a concrete example, for the densities of
$N(\theta,1)$, $\theta\in\Reals$, $A(\theta, \theta', x)\equiv 1$ and
hence Assumption 4 is satisfied.

As $\theta-\theta_0\to 0$, the nulls become increasingly similar.  By
the asymptotic theory of statistics, to attain a fixed power while
keeping the same significance level for the tests, $k$ should grow at
the same order as $(\theta-\theta_0)^{-2}$.  On the other hand, the
results below deal with the case where $k$ grows a little more slowly.
First consider the univariate case $d=1$.
\begin{theorem} \label{thm:lrt-1d}
  Let $Q_\alpha>1$.  Denote $\delta = \theta-\theta_0$ and $k$ the
  number of replications per null.  Suppose 
  \begin{align}
    k = \frac{1}{\delta^2 s(\delta)}, \ \text{ such that }\ 
    s(\delta)\toi \ \text{ and }\
    \frac{s(\delta)}{\ln (1/\delta^2)} \to 0 \ \text{ as }\ \delta\to
    0.
    \label{eq:lrt-k}
  \end{align}
  Then, as $\theta\to\theta_0$, 
  \begin{gather}
    p_{k,\delta}(\alpha)=
    (1+o(1)) 
    \sqrt{\frac{(1-\pa)\pa}{2\pi (1-\alpha) \alpha}}
    \frac{\sqrt{k I(\theta_0)} \delta}{\ln Q_\alpha}
    \exp\Cbr{
      -\frac{(\ln Q_\alpha)^2}{2 k I(\theta_0) \delta^2}
    }.
    \label{eq:p-lrt}
  \end{gather}
\end{theorem}

Note that $Q_\alpha>1$ if $\pa + \alpha < 1$.  In practice, since
$\pa$ is usually much less than 1 and $\alpha$ is small or only
moderately large, the assumption $Q_\alpha>1$ is not restrictive.

The multivariate case $d>1$ can be derived as a corollary.
\begin{cor} \label{cor:lrt}
  Let $Q_\alpha>1$.  Denote $\delta = \theta - \theta_0$.  Suppose $k$
  satisfies \eqref{eq:lrt-k}, with $\delta^2$ being replaced by
  $|\delta|^2$.  Let $q(\delta) = \delta^T I(\theta_0) \delta$.
  Then, as $\theta\to\theta_0$,
  \begin{gather}
    p_{k,\delta}(\alpha) = (1+o(1))
    \sqrt{\frac{(1-\pa) \pa}{2\pi (1-\alpha)\alpha}}
    \frac{\sqrt{k q(\delta)}}{\ln Q_\alpha}
    \exp\Cbr{-
      \frac{(\ln Q_\alpha)^2}{2 k q(\delta)}
    }.
    \label{eq:p-lrt-mv}
  \end{gather}
\end{cor}

Under the condition of Theorem \ref{thm:lrt-1d}, it is not hard to see
$p_{k,\delta}(\alpha)\to 0$.  Therefore, the minimum number of nulls
$N_*$ and the minimum volume of data $V_*$ in order for
\eqref{eq:detect} to be satisfied have the asympotics in
\eqref{eq:min-n}, yielding
$$
V_* = k N_* = (1+o(1)) \ln\frac{1}{1-p} \sqrt{\frac{2\pi (1-\alpha)
    \alpha}{(1-\pa)\pa}} \frac{ \sqrt{k}\ln
  Q_\alpha}{\sqrt{I(\theta_0)} \delta} \exp\Cbr{
  \frac{(\ln Q_\alpha)^2}{2 k I(\theta_0) \delta^2}
}.
$$

On the other hand, if $k$ has the same order as $\delta^{-2}$, then by
the Central Limit Theorem, in order to satisfy \eqref{eq:detect},
$N_*$ only needs to be a large constant and $V_*$ is of the same
order as $\delta^{-2}$.  To see in which case the minimum data volume
is larger, it suffices to compare the orders of
$\frac{\sqrt{k}}{\delta} \exp\Cbr{\frac{c}{k\delta^2}}$ and
$\delta^{-2}$ as $\delta\to 0$, where $c>0$ is a constant.  By
$k=\frac{1}{\delta^2 s(\delta)}$ with $s(\delta)\toi$, the ratio of
the two is $\delta\sqrt{k} \exp\{\frac{c}{k\delta^2}\}
=e^{c s(\delta)}/\sqrt{s(\delta)} \toi$.  Therefore, when $k$ cannot
grow as fast as $\delta^{-2}$, a much larger volume of data is
required.

\subsection{Multiple tests on means of normal distributions}
\label{subsec:normal}
Consider nulls $H_i: \theta_i=\theta_0$ for $N(\theta_i, \sigma^2)$,
where $\sigma^2$ is known and under false $H_i$, $\theta_i=\theta$,
with $\theta-\theta_0=\delta>0$.  Without loss of generality, let
$\theta_0=0$ and $\sigma=1$.  Then $f_0$ and $f_a$ are the densities
of $N(0,1)$ and $N(\delta, 1)$, respectively.  By $f_a(x)/f_0(x) =
\exp(\delta x - \nth{2} \delta^2)$,  the event in \eqref{eq:pfdr}
becomes
$$
\cE_i(k) =  \Cbr{\sum_{j=1}^k X_{ij} \ge \frac{\ln Q_\alpha}{\delta} +
  \frac{k\delta}{2}
}.
$$

Under true $H_i$, $\sum_j X_{ij} \sim \sqrt{k} Z$, while under false
$H_i$, $\sum_j X_{ij} \sim \sqrt{k} Z+k\delta$, where $Z\sim N(0,1)$.
Therefore, by \eqref{eq:pfdr2},
\begin{align}
  p_{k,\delta}(\alpha) = (1-a) \bar\Phi\Grp{
    \frac{\ln Q_\alpha}{\sqrt{k} \delta} + \frac{\sqrt{k}\delta}{2}
  }+
  a \bar\Phi\Grp{
    \frac{\ln Q_\alpha}{\sqrt{k} \delta} - \frac{\sqrt{k}\delta}{2}
  },  \label{eq:tail-normal}
\end{align}
where $\bar\Phi(t) = 1-\Phi(t) = 1-P(Z\le t)$.

Let $Q_\alpha>1$.  If $\delta\downarrow 0$ such that
$k\delta^2\to 0$, then $\frac{\ln Q_\alpha}{\sqrt{k} \delta}\pm
\frac{\sqrt{k}\delta}{2} \to\infty$.  Recall that, as $t\toi$,
\begin{align}
  \bar\Phi(t) = \Phi(-t) = (1+o(1))
  \frac{e^{-t^2/2}}{\sqrt{2\pi}\,t}.
  \label{eq:normal}
\end{align}
It is then not hard to get the asymptotic in Theorem \ref{thm:lrt-1d}
with fewer restrictions on $(\delta, k)$.
\begin{cor} \label{cor:normal}
  Let $Q_\alpha>1$.  Suppose $k\delta^2\to 0$ as $\delta\to 0$.  Let
  $f_\theta$ be the densities of $N(\theta, \sigma^2)$, with
  $\sigma^2$ being known.  Then, as $\theta\to\theta_0$,
  \begin{gather}
    p_{k,\delta}(\alpha)\sim
    \sqrt{\frac{(1-\pa)\pa}{2\pi (1-\alpha) \alpha}}
    \frac{\sqrt{k} \delta/\sigma}{\ln Q_\alpha}
    \exp\Cbr{
      -\frac{(\ln Q_\alpha)^2}{2 k\delta^2/\sigma^2}
    }, \label{eq:norm-P}\\
    V_* = kN_* \sim \ln\nth{1-p}\sqrt{
      \frac{2\pi(1-\alpha)\alpha}{(1-\pa)\pa}}
    \frac{\sqrt{k} \ln Q_\alpha}{\delta/\sigma} \exp
    \Cbr{\frac{(\ln Q_\alpha)^2}{2k\delta^2/\sigma^2}}.  \label{eq:norm-V}
  \end{gather}
\end{cor}

\begin{figure}[t]
  \renewcommand{\arraystretch}{0}
  \begin{center}
    \setlength{\unitlength}{1mm}
    \begin{picture}(155,50)(-5,5)
      \put(0,0){\epsfig{file = 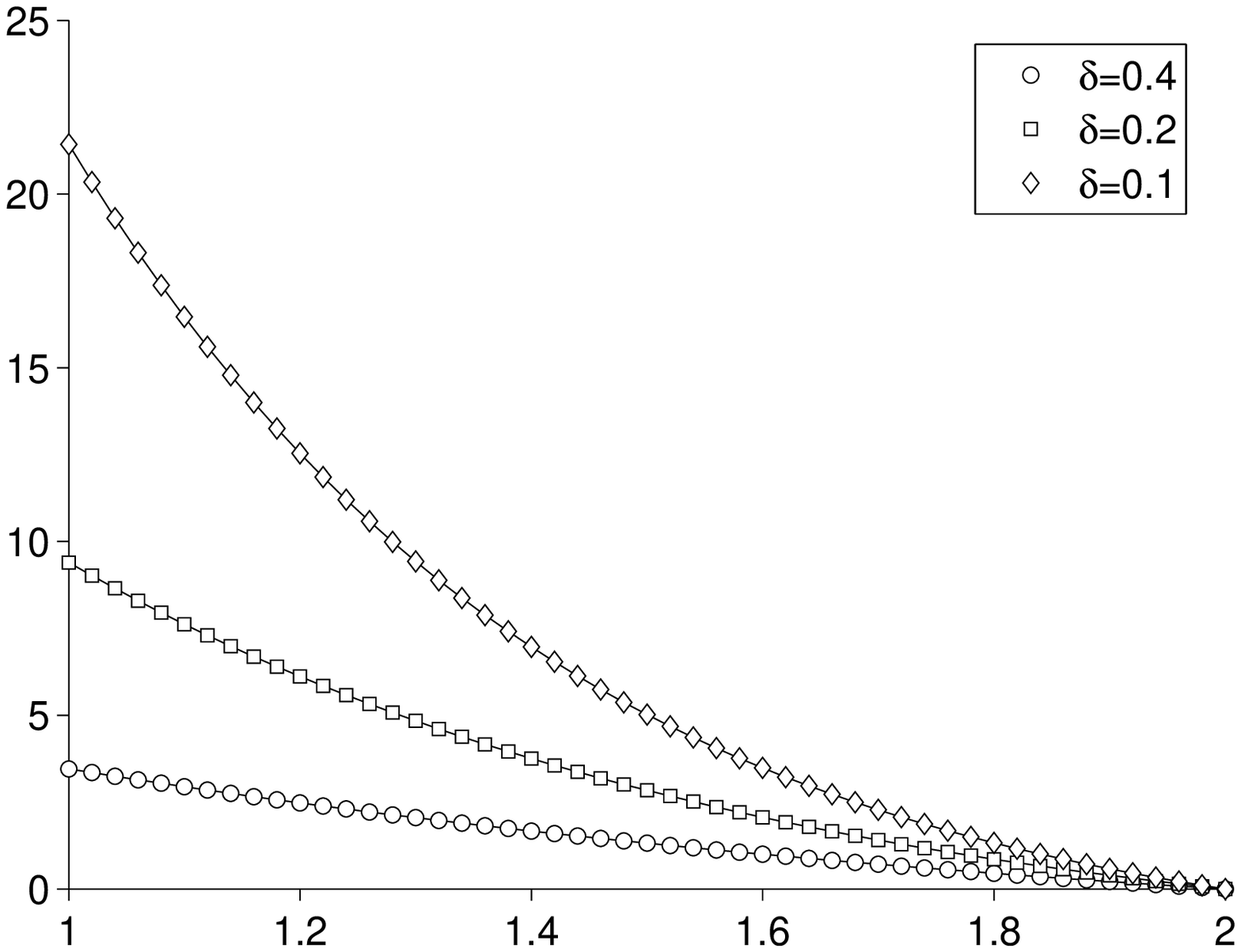, width=7cm}}
      \put(-5,17){\rotatebox{90}{\small $\log_{10}(V_t/V_2)$}}
      \put(80,0){\epsfig{file = 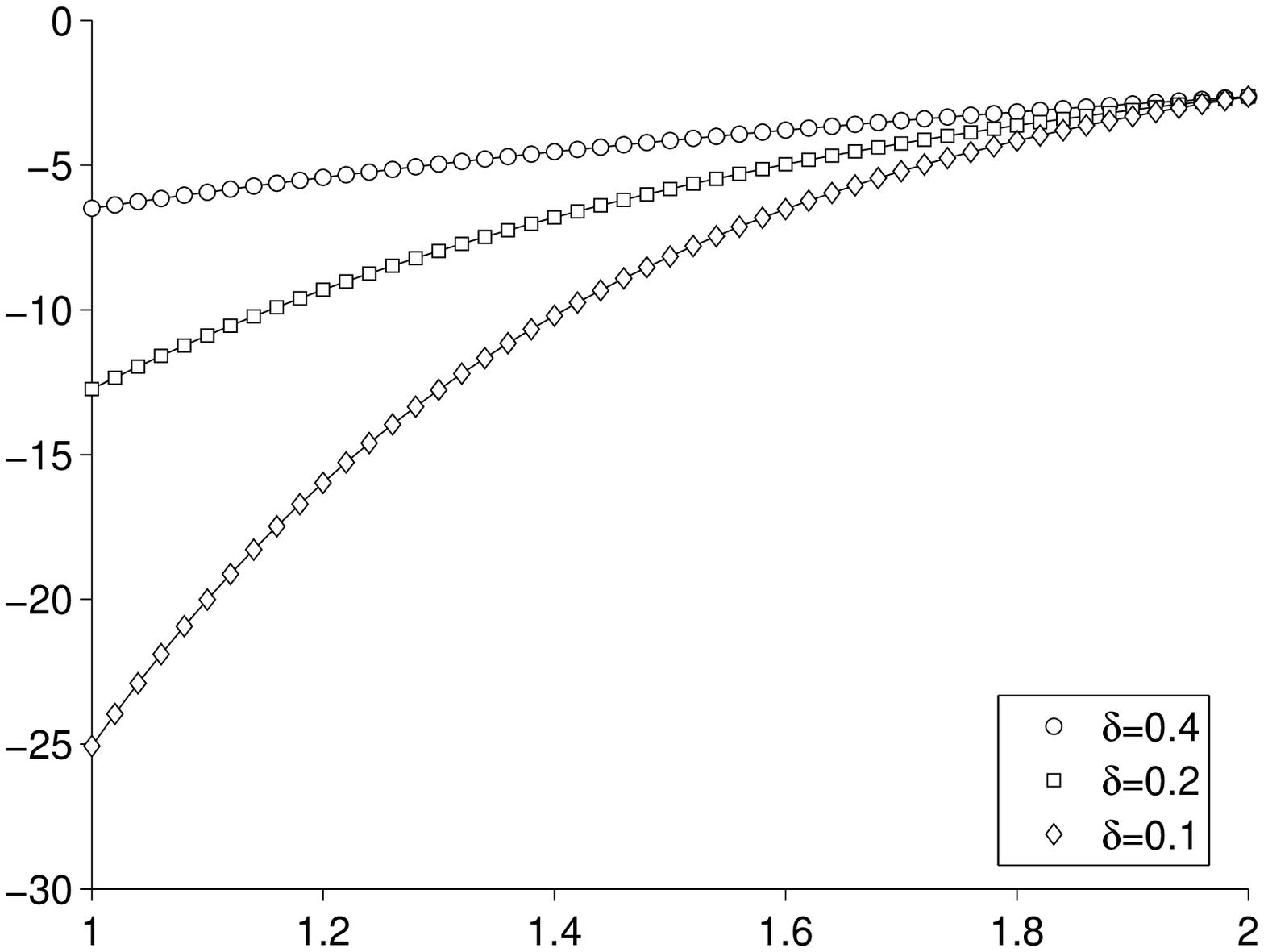, width=7cm}}
      \put(75,21){\rotatebox{90}{\small $\log_{10} P_t$}}
      \put(36,-2.5){\small $t$} \put(117,-2.5){\small $t$}
      \put(35,50){\small (A)} \put(116,50){\small (B)}
    \end{picture}
  \end{center}
  \begin{small}
    \caption{Minimum data volume and power for tests on mean
      values of $N(\delta,1)$ when $k=\delta^{-t}$, $t\in [1,2]$; see
      details in Sections \ref{subsec:normal} and \ref{subsec:power}.}
  \end{small}
  \label{fig}
\end{figure}

The rapid increase of minimum data volume is illustrated in Figure
\ref{fig}(A), which graphs $\log_{10}(V_t/V_2)$ versus
$t\in [0,1]$ for $\delta=0.1$, $0.2$, and $0.4$, where $V_t$ is the
right hand side of \eqref{eq:norm-V} with $k=\delta^{-t}$.  For the
plot, $\pa = 5\%$, $\alpha = 0.4$ and $p=0.9$.  Even at the log scale, 
the increase in the minimum data value is apparent.

\subsection{Multiple tests on scales of Gamma distributions}
\label{subsec:gamma}
Denote by $\dgamma(a,b)$ the Gamma distribution with shape parameter
$a$ and scale parameter $b$.  Multiple tests on the scales of Gamma
distributions have been used as a case of study in the literature
\citep{donoho:jin:04}.  Fix $\shape>0$.  Let $f_0(x)$ be the density of
$\dgamma(\shape,1)$ and $f_a(x)$ that of $\dgamma(\shape, 1+\delta)$,
where $\delta>0$.  Then
\begin{align*}
  f_0(x) = \frac{x^{\shape-1} e^{-x}}{\Gamma(\shape)}, \quad
  f_a(x)=f_\delta(x) = \frac{x^{\shape-1}e^{-x/(1+\delta)}}
  {\Gamma(\shape) (1+\delta)^\shape} , \quad x>0.
\end{align*}
By $f_a(x)/f_0(x) = (1+\delta)^{-\shape} e^{\delta x/(1+\delta)}$, the
event in \eqref{eq:pfdr} becomes
\begin{align}
  \cE_i(k)&
  = 
  \Cbr{\sum_{j=1}^k X_{ij} \ge c_k},\
  \text{ with }
  c_k:= \frac{[k\shape \ln(1+\delta) + \ln
    Q_\alpha](1+\delta)}{\delta}.
  \label{eq:gamma-thresh}
\end{align}

Under true $H_i$, $\sum_j X_{ij} \sim \dgamma(k\shape,1)$;
under false $H_i$, $\sum_j X_{ij} \sim \dgamma(k\shape,1+\delta)$.  By
\eqref{eq:pfdr2}, to get the asymptotics of $N_*$ and $V_*$, the main
step is to get the asymptotics of the probability of $\{S \ge c_k\}$
for $S$ following $\dgamma(k\shape,1)$ or $\dgamma(k\shape,1+\delta)$.
Since the tail probabilities under Gamma distributions are available
in detail, the asymptotics can be attained for a much wider range of
values of $k$ than for the general case.  The results are as follows;
see Section \ref{sec:proof-gamma} for a proof.
\begin{theorem} \label{thm:lrt-gamma}
  Let $f_0$ and $f_a$ be as above and $Q_\alpha>1$.  Suppose
  \begin{align}
    (\delta, k)\to (0,\infty)\ \text{ such that }\
    k\delta \toi, \ k\delta^2\to 0.
    \label{eq:k-delta}
  \end{align}
  Then, denoting $\psi(t) = t-\frac{t^2}{2}-\ln(1+t)$ for $t>-1$,
  \begin{gather}
    p_{k,\delta}(\alpha)=(1+o(1))
    \sqrt{\frac{(1-\pa)\pa}{2\pi(1-\alpha)\alpha}}
    \frac{\sqrt{k\shape} \delta}{\sqrt{\ln Q_\alpha}}
    \exp\Cbr{
      -\frac{(\ln Q_\alpha) ^2}{2k\shape \delta^2}
      - k\shape \psi\Grp{
        \frac{\ln Q_\alpha}{k\shape\delta}}
    }.
    \label{eq:p-gamma}
  \end{gather}
\end{theorem}

\subsection{Asymptotics of power} \label{subsec:power}
For fixed $\delta$ and $k$, power can be analyzed
using previous results \citep{signorovitch:06, storey:07:opt, chi:08}.
Under the setup here, since $(\delta, k)\to (0,\infty)$, the
asymptotics of power are of interest.  To avoid subtleties that a
finite number of nulls may cause, we consider power 
under the situation where arbitrarily many nulls can be
tested.  For any procedure, let $N_a$ and $R_a = R-R_0$ denote the
numbers of false nulls and rejected false nulls, respectively.
Then, provided the limit below exists,
\begin{align}
  \power_\infty = \lim_{N\toi} E[R_a/N_a],  \label{eq:power}
\end{align}
characterizes the power of the procedure when $N\gg N_*$.  We compare
the powers of different procedures when they control the pFDR around
or below the same level.  For $N\gg N_*$, the pFDR of a procedure 
can be characterized by 
\begin{align}
  \pfdr_\infty = \lim_{N\toi} E[R_0/R\gv R>0].  \label{eq:pfdr-infty}
\end{align}

As the limits in \eqref{eq:power} and \eqref{eq:pfdr-infty} are
defined for fixed $\delta$ and $k$, we use $\power_\infty(\delta, k)$
and $\pfdr_\infty(\delta, k)$ to indicate the dependency and consider
the asymptotics of the quantities as 
$(\delta,k)\to (0,\infty)$. \comment{In view of the criterion
\eqref{eq:criterion}, another possible way to define 
$\power_\infty$ and $\pfdr_\infty$ is to use conditional expectations,
e.g., $E[R_a/N_a\gv \vf X]$ in place of $E[R_a/N_a]$.  Under the
random effects model \eqref{eq:model}, for regular procedures, the
conditional expectations converge in probability to the same limits of
their unconditional counterparts, so using conditional expectations
makes no difference to the discussion.}
Fix $\alpha\in (0,1)$.  First consider the thresholding procedure with
cut-off $\alpha$,
\begin{align*}
  d_i^*(\vf X) = \cf{\Pr(\eta_i=0\gv\vf X)\le\alpha}, \quad
  i=1,2,\ldots, 
\end{align*}
i.e., $d^*$ rejects $H_i$ if and only $\Pr(\eta_i=0\gv \vf X)\le
\alpha$.  Denote by $\power_\infty^*(\delta, k)$ and
$\pfdr_\infty^*(\delta, k)$ the asymptotic power and pFDR of $d^*$,
respectively.

\begin{prop} \label{prop:odds-ratio}
  Suppose $(\delta, k)\to (0,\infty)$ such that 
  \begin{align}
    p_{k,\delta}(\alpha) \to 0\ \text{ while staying positive, and }
    \label{eq:positive} \\
    p_{k,\delta}(\alpha_1) = o(p_{k,\delta}(\alpha)) \ \text{ for any
    }\ 0<\alpha_1<\alpha,
    \label{eq:increase}
  \end{align}
  Then $d^*$ has the following property
  \begin{align}
    \begin{array}{c}
      \text{for any fixed $\delta>0$ and $k\ge 1$, the limits
        in \eqref{eq:power} and \eqref{eq:pfdr-infty} exist,} \\
      \text{and $R_a/N_a\convP \power_\infty$ as $N\toi$.}
    \end{array}
    \label{eq:property}
  \end{align}
  Moreover,
  \begin{align}
    \power_\infty^*(\delta,k) = (1+o(1))
    \frac{(1-\alpha)p_{k,\delta}(\alpha)}{\pa}, \quad
    \pfdr_\infty^*(\delta,k) = (1+o(1)) \alpha.
    \label{eq:odds-large-n}
  \end{align}
\end{prop}

We use $d^*$ as a baseline to compare other procedures that satisfy
the basic property \eqref{eq:property} while asymptotically
controlling the pFDR.
\begin{prop}\label{prop:large-n}
  Let $(\delta, k)\to(0,\infty)$ as in Proposition
  \ref{prop:odds-ratio}.  Let $d$ be a procedure satisfying
  \eqref{eq:property} with $\Limsup \pfdr_\infty(\delta,k)\le
  \alpha$.  If $\,\power_\infty(\delta,k) \ge \power_\infty^*(\delta,
  k)$ for all $(\delta,k)$, then for any $\alpha_2>\alpha$,
  \begin{align}
    \power_\infty(\delta,k) \le
    (1+o(1))\frac{p_{k,\delta}(\alpha_2)}{\pa},
    \quad
    \pfdr_\infty(\delta,k)=(1+o(1))\alpha.
    \label{eq:large-n}
  \end{align}
\end{prop}

It is not hard to see that provided $Q_\alpha>1$, the
$p_{k,\delta}(\alpha)$ given in Theorems \ref{thm:lrt-1d} and
\ref{thm:lrt-gamma} satisfies \eqref{eq:positive} and
\eqref{eq:increase} and therefore the above results apply.  Since by
Proposition \ref{prop:odds-ratio}, $p_{k,\delta}(\alpha_2)$ is of the
same order as the power of a thresholding procedure with cut-off
$\alpha_2$, Proposition \ref{prop:large-n} basically says that for any
procedure satisfying \eqref{eq:property} with $\Limsup
\pfdr_\infty(\delta,k)\le\alpha$, its power is dominated up to a
constant factor by a thresholding procedure with a cut-off just a
little bit above $\alpha$.  In view of this, one question is whether
there is a most powerful procedure among those that satisfy
\eqref{eq:property} and $\Limsup\pfdr_\infty(\delta,k) \le \alpha$,
and in particular, whether $d^*$ is such one.  As seen next, in
general the answer is no.  Given $c>0$, let $d$ be a procedure such
that for each $(\delta, k)$, its cut-off is $\alpha+ck\delta^2$, i.e.
$$
d_i(\vf X) = \cf{\Pr(\eta_i=0\gv \vf X)\le \alpha + c k\delta^2}.
$$

\begin{prop} \label{prop:power}
  Under the random effects model \eqref{eq:model}, let $(\delta,k)$ be
  as in Theorem \ref{thm:lrt-1d}.  Suppose $Q_\alpha>1$.  Given $M>1$,
  let $c = (1-\alpha)\alpha I(\theta_0) \ln M/\ln Q_\alpha$.  Then $d$
  satisfies \eqref{eq:property} and $\Limsup\pfdr_\infty(\delta,k)\le
  \alpha$, while $\power_\infty(\delta, k) = (M+o(1))
  \power^*_\infty(\delta,k)$.

  More generally, there are no asymptotically most powerful procedures
  that satisfy \eqref{eq:property} and
  $\Limsup\pfdr_\infty(\delta,k)\le \alpha$.
\end{prop}

As an illustration, consider multiple testing for the mean values of
$N(\theta,1)$ as in Section \ref{subsec:normal}.  Figure~\ref{fig}(B)
shows the dependency of the asymptotic power of $d^*$ on $(\delta,k)$.
Under the same parameters as in panel (A), it graphs $\log_{10}P_t$,
$t\in [1,2]$, where $P_t$ is $(1-\alpha)/\alpha$ times the right hand
side of \eqref{eq:norm-P} with $k=\delta^{-t}$.  From Proposition
\ref{prop:odds-ratio}, we know $\power_\infty^*(\delta, k) = (1+o(1))
P_t$ as $\delta\to 0$ and $k=\delta^{-t}$.  The rapid decrease of
power as $\delta\to 0$ is clear.  We next illustrate how the
asymptotic power of a thresholding procedure can be arbitrarily
increased by a small change in cut-off.  As seen from
\eqref{eq:tail-normal}, for any $\delta$ and $k$,
$\power_\infty^*(\delta,k)=\bar\Phi\Grp{\frac{\ln Q_\alpha}{\sqrt{k}
    \delta} - \frac{\sqrt{k}\delta}{2}}$ and the thresholding
procedure with cut-off $\alpha+ck\delta^2$ has
$\power_\infty(\delta,k) = \bar\Phi\Grp{\frac{\ln Q_{\alpha +
      ck\delta^2}}{\sqrt{k} \delta} - \frac{\sqrt{k}\delta}{2}}$.  As
$(\delta, k)\to (0,\infty)$ with $k\delta^2\to 0$, the difference
between the cut-offs $\alpha$ and $\alpha+ck\delta^2$ tends to 0.  It
is not hard to get that for both procedures, $\pfdr_\infty(\delta,k)
\to\alpha$.  On the other hand, by \eqref{eq:normal}, the ratio of
their asymptotic powers is
\begin{align*}
  &
  (1+o(1))\exp\Cbr{
    \Grp{
      \frac{\ln Q_{\alpha + ck\delta^2}}{
        \sqrt{k} \delta} - \frac{\sqrt{k}\delta}{2}
    }^2 -
    \Grp{
      \frac{\ln Q_\alpha}{
        \sqrt{k} \delta} -\frac{\sqrt{k}\delta}{2}
    }^2
  } \\
  =
  &
  (1+o(1)) \exp\Cbr{\frac{(\ln Q_{\alpha+ck\delta^2})^2 -
      (\ln Q_\alpha)^2}{k\delta^2}}
  =(1+o(1)) \exp\Cbr{-\frac{2c\ln Q_\alpha}{1-\alpha}}.
\end{align*}
Therefore, by increasing $c$, the power of the second thresholding
procedure is arbitrarily many times higher than $d^*$.

\section{Summary and remarks} \label{sec:summary}
This article studies the issues of minimum data volume and power when
$k= o(\delta^{-2})$, i.e., the number of repeated measurements for
each null is much smaller than the squared differences between false
and true nulls.  It shows that in this case, in order to meet a pFDR
based performance criterion, the minimum data volume has to grow much
faster than in the case where $k$ is of the same order as
$\delta^{-2}$.  It also shows how fast power will decay to 0 and the
sensitivity of the power to small changes in rejection rules.

The results are essentially due to the fact that when $k$ is not large
enough, evidence against true nulls can only come from values of test
statistics far away from the ``normal'' ones.  When $k$ increases more
slowly than $\delta^{-2}$ but faster than $\delta^{-1}$, such values
can be treated as moderate deviations \citep{dembo:zeitouni:ldt:98},
which can yield the log-growth rate of the minimum data volume but
nevertheless may not be accurate enough to give the growth rate
itself.  On the other hand, the article does not consider the case
where $k$ is only of the order of $\delta^{-1}$.  Analysis in this
case seems to require the large deviations principle and can be quite
subtle \citep{chi:07b,dembo:zeitouni:ldt:98}.

\bibliographystyle{agsm}

\renewcommand{\thesection}{A}
\section*{Appendix: technical details} \label{sec:proofs}
\subsection{Proof for the equivalence of criteria \eqref{eq:criterion}
  and \eqref{eq:detect}}
We sketch a proof that under the random effects model
\eqref{eq:model}, the criterion \eqref{eq:criterion} can be rewritten
as \eqref{eq:detect}, where the infimum in \eqref{eq:criterion} is
taken over procedures satisfying \eqref{eq:multitest}.  For more
details, see \cite{chi:tan:08}.

For any procedure $d(\vf X)$ as in \eqref{eq:multitest},
$R=\sum_{i=1}^N d_i(\vf X)$ and $R_0 = \sum_{i=1}^N \cf{\eta_i=0}
d_i(\vf X)$.  Given $\vf X$, if $R>0$, then, as $d_i(\vf X)$ are now
determined,
\begin{align*}
  E\Sbr{R_0/R\gv \vf X}
  = \frac{1}{R} \sum_{i=1}^N d_i(\vf X) \Pr(\eta_i = 0 \gv \vf X)
  \ge \min_{i=1}^N \Pr(\eta_i=0\gv\vf X),
\end{align*}
with equality if and only if $d$ only rejects $H_i$ with the smallest
$\Pr(\eta_i = 0\gv \vf X)$.  On the other hand, if $R=0$, then by
definition, $E[R_0/R\gv R>0, \vf X]=1$.  Note that by Bayes rule, 
\begin{align}
  \Pr(\eta_i = 0 \gv \vf X)
   = \Pr(\eta_i = 0\gv X_{i j},\, j=1,\ldots, k) 
  = \Sbr{1+\frac{\pa}{1-\pa} \prod_{j=1}^k
    \frac{f_a(X_{i j})}{f_0(X_{i j})}
  }^{-1}.  \label{eq:posterior}
\end{align}
It is then seen that the criterion \eqref{eq:criterion} can be
rewritten as \eqref{eq:detect}.

\subsection{Proofs for general multiple tests}
Recall that by Bikjalis' theorem \citep{nagaev:79}, there is an
absolute constant $\beta>0$, such that for any $\inum Z$ i.i.d.\
with $EZ_1=0$, $\var(Z_1)=\sigma^2>0$ and
$E|Z_1|^3<\infty$,
\begin{align}
  \Abs{
    \bar\Phi(t) - \Pr\Grp{
      \nth{\sigma\sqrt{k}} \sum_{i=1}^k Z_i \ge t
    }
  } \le \frac{\beta E|Z_1|^3}{\sigma^3 \sqrt{k}(1+|t|^3)}, \quad
  k=1,2,\ldots.
  \label{eq:mdp}
\end{align}
\paragraph{\it Proof of Theorem \ref{thm:lrt-1d}.\ }
  For $\theta\in \Theta\subset \Reals$, let $\delta =
  \theta-\theta_0$.  Then
\begin{align}
  \ell(\theta, x) - \ell(\theta_0, x)&
  =
  \dot\ell(\theta_0, x)\delta
  + \frac{A(\theta_0, \theta, x) \delta^2}{2} \label{eq:taylor0} \\
  &
  =
  \dot\ell(\theta, x)\delta
  - \frac{A(\theta, \theta_0, x) \delta^2}{2}.  \label{eq:taylor}
\end{align}

According to \eqref{eq:min-n}, we need to compute
$p_{k,\delta}(\alpha)$.  By \eqref{eq:pfdr2},
\begin{align}
  p_{k,\delta}(\alpha) = (1-\pa) P_{\theta_0}(E_k) + \pa P_\theta(E_k),
  \label{eq:p-mix}
\end{align}
where $P_\theta$ is the $k$-fold product of the probability measure with
density $f_\theta$ and $E_k$ is the event $\{(\eno X k): \sum_{i=1}^k
[\ell(\theta, X_i) - \ell(\theta_0, X_i)] \ge \ln Q_\alpha\}$.  Under
$f_{\theta_0}$, $\dot\ell(\theta_0, X_i)$ are i.i.d.\  with mean
0 and variance $I(\theta_0)$.  Given $\epsilon>0$, define events
\begin{align*}
  G_k =\Cbr{
    \Abs{
      \nth k\sum_{i=1}^k A(\theta_0, \theta, X_i) + I(\theta_0)
    } \le \epsilon
  }.
\end{align*}
Denote $Z_i = \dot\ell(\theta_0, X_i)$.  By \eqref{eq:taylor0},
$F_-\cap G_k \subset   E_k\cap G_k  \subset F_+\cap G_k$, where, for
$\theta\not=\theta_0$,
\begin{align*}
  F_\pm
  =\Cbr{
    \frac{\text{sign}(\delta)}{\sqrt{k I(\theta_0)}}
    \sum_{i=1}^k Z_i
    \ge u_\pm
  } \
  \text{ with }\
  u_\pm
  = \frac{
    \ln Q_\alpha + \nth 2 k(I(\theta_0)\mp\epsilon)\delta^2
  }{
    \sqrt{k I(\theta_0)}\, |\delta|
  }.
\end{align*}
Without loss of generality, we only consider the case
$\theta>\theta_0$.  Then
\begin{align}
  P_{\theta_0}(F_-) - P_{\theta_0}(G_k^c) \le P_{\theta_0}(E_k) \le
  P_{\theta_0}(F_+) + P_{\theta_0}(G_k^c).
  \label{eq:prob}
\end{align}
By \eqref{eq:mdp} and Assumption \ref{A:mnt},
\begin{align*}
  \Abs{
    \bar\Phi(u_\pm) - P_{\theta_0}(F_\pm)
  }
  &
  \le \frac{\beta M_3 }
  {\sqrt{k}\,I(\theta_0)^{3/2} (1+|u_\pm|^3)}.
\end{align*}
We need the following results.
\begin{lemma} \label{lemma:lrt}
  If $k$ satisfies \eqref{eq:lrt-k}, then, as $\delta\to 0$,
  \begin{gather}
    \bar\Phi(u_\pm) = (1+o(1)) \frac{e^{-u^2_\pm/2}}{\sqrt{2\pi}
      u_\pm}, \label{eq:phi} \\
    \nth{
      \sqrt{k} (1+|u_\pm|^3)
    } = o(\bar\Phi(u_\pm)), \label{eq:approx1} \\
    P_{\theta_0}(G_k^c) = o(\bar\Phi(u_\pm)).
    \label{eq:approx2}
  \end{gather}
\end{lemma}

Assume Lemma \ref{lemma:lrt} is true for now.  By
\eqref{eq:prob}--\eqref{eq:approx2}, there is $r_k\to 0$, such that
\begin{align*}
  (1-r_k) \frac{e^{-u_-^2/2}}{\sqrt{2\pi} u_-}
  \le P_{\theta_0}(E_k) \le 
  (1+r_k) \frac{e^{-u_+^2/2}}{\sqrt{2\pi} u_+}.
\end{align*}
Since $k\delta^2\to 0$, $u_-\sim u_+ \sim \frac{\ln Q_\alpha}{
    \sqrt{k I(\theta_0)}\,\delta}$.  On the other hand,
\begin{align*}
  u_\pm^2 =
  \frac{(\ln Q_\alpha)^2}{k I(\theta_0) \delta^2}
  +\frac{(I(\theta_0)\mp \epsilon) \ln Q_\alpha}{I(\theta_0)}
  + \frac{
    (I(\theta_0) \mp \epsilon)^2 k \delta^2}{4 I(\theta_0)}.
\end{align*}
By Assumption \ref{A:I-pd}, $I(\theta_0)>0$.  Since $\epsilon>0$ is
arbitrary and $k\delta^2\to 0$, we then get
\begin{align}
  P_{\theta_0}(E_k) = (1+o(1)) 
  \frac{\sqrt{k I(\theta_0)}\,\delta}{\sqrt{2\pi}
    \ln Q_\alpha
  }
  \exp\Cbr{
    -\frac{(\ln Q_\alpha)^2}{2 k I(\theta_0) \delta^2}
    - \frac{\ln Q_\alpha}{2}
  }.
  \label{eq:p-lrt0}
\end{align}

With similar argument, now applied to $\dot\ell(\theta, X_i)$ under
$f_\theta$,
\begin{align*}
  P_\theta(E_k) = (1+o(1)) 
  \frac{\sqrt{k I(\theta)}\delta}{\sqrt{2\pi} \ln Q_\alpha}
  \exp\Cbr{
    -\frac{(\ln Q_\alpha)^2}{2 k I(\theta) \delta^2}
    + \frac{\ln Q_\alpha}{2}
  },
\end{align*}
where $+\nth 2\ln Q_\alpha$ in the exponential is
due to $-\nth 2 A(\theta, \theta_0, x)\delta^2$ in
\eqref{eq:taylor}.  By Assumption \ref{A:I-diff}, there are constants
$C>0$ and $r>0$, such that for $|\delta|<r$,
\begin{align*}
  \Abs{
    \frac{1}{k I(\theta) \delta^2} - \frac{1}{k I(\theta_0) \delta^2} 
  }
  \le \frac{C}{ k|\delta| I^2(\theta_0)}.
\end{align*}
Since $k\delta\toi$, it follows that
\begin{align}
  P_\theta(E_k) = 
  (1+o(1)) 
  \frac{\sqrt{ k I(\theta_0)} \delta}{\sqrt{2\pi}
    \ln Q_\alpha}
  \exp\Cbr{
    -\frac{(\ln Q_\alpha)^2}{2 k I(\theta_0) \delta^2}
    + \frac{\ln Q_\alpha}{2}
  }.
  \label{eq:p-lrt1}
\end{align}
Combining \eqref{eq:p-mix}, \eqref{eq:p-lrt0}, \eqref{eq:p-lrt1}, and
$\exp\{\frac{\ln Q_\alpha}{2}\} =
\sqrt{\frac{(1-\pa)(1-\alpha)}{\pa\alpha}}$, 
\eqref{eq:p-lrt} then follows.  \qed

\paragraph{\it Proof of Lemma \ref{lemma:lrt}.\/}  Because $u_\pm\toi$ as
$\delta\to 0$, \eqref{eq:phi} follows from \eqref{eq:normal}.  To show
\eqref{eq:approx1}, it suffices to show $\sqrt{k} u_\pm^2
e^{-u_\pm^2/2} \toi$ as $\delta\to 0$, or, equivalently,
$\ln k - u^2_\pm + 4\ln u_\pm \toi$.  Because $u_\pm$ is of the same
order as $\frac{1}{\sqrt{k}\delta}$ and 
$k\delta^2\to 0$, it is seen the above asymptotic follows if
$\frac{1}{k\delta^2}= o(\ln k)$, or $(k\ln k)\delta^2\toi$.  Now by
$s(\delta)\toi$ and $s(\delta) = o(\ln (1/\delta))$,
$$
(k\ln k)\delta^2 = \frac{1}{\delta^2 s(\delta)}
\Grp{2\ln\frac{1}{\delta} - \ln s(\delta)} \delta^2 \toi.
$$
To show \eqref{eq:approx2}, let $\lambda>0$ be as in \eqref{eq:LDP}.
Then $P_{\theta_0}(G_k^c)$ is of the same order as $e^{-\lambda k}$.
By \eqref{eq:phi}, it suffices to show $u^2_\pm = o(k)$.  Since
$u^2_\pm$ is of the same order as $1/(k\delta^2)$ and $k\delta\toi$,
the last claim is proved.  \qed

\paragraph{\it Proof of Corollary \ref{cor:lrt}.\/}  In place of
\eqref{eq:taylor0} and \eqref{eq:taylor}, we have
\begin{align*}
  \ell(\theta, x) - \ell(\theta_0, x)
  =
  \dot\ell(\theta_0, x)^T\delta
  + \frac{\delta^T A(\theta_0, \theta, x) \delta}{2} 
  =
  \dot\ell(\theta, x)^T\delta
  - \frac{\delta^T A(\theta, \theta_0, x) \delta}{2}.
\end{align*}
Let $e = \delta / |\delta|$.  Under $f_{\theta_0}$,
$\dot\ell(\theta_0, X_i)^T e$ are i.i.d.\ with mean 0 and variance
$v(\theta_0)=e^T I(\theta_0) e>0$; under $f_\theta$,
$\dot\ell(\theta, X_i)^T e$ are i.i.d.\ with mean 0 and variance
$v(\theta) = e^T I(\theta) e$.  Applying the proof of Theorem
\ref{thm:lrt-1d} to $\dot\ell(\theta_0, X_i)^T e$ and
$\dot\ell(\theta, X_i)^T e$ yields
\begin{align*}
  p_{k,\delta}(\alpha) \sim 
  \sqrt{\frac{(1-\pa)\pa}{2\pi (1-\alpha) \alpha}}
  \frac{\sqrt{k v(\theta_0)}\,|\delta|}{\ln Q_\alpha}
  \exp\Cbr{
    -\frac{(\ln Q_\alpha)^2}{2 k v(\theta_0)\,|\delta|^2}
  }.
\end{align*}
Since $v(\theta_0)|\delta|^2 = \delta^T I(\theta_0)\delta$,
\eqref{eq:p-lrt-mv} then follows.  \qed

\subsection{Proofs for multiple tests on Gamma distributions}
\label{sec:proof-gamma} 
We next prove Theorem \ref{thm:lrt-gamma}.  Denote by $G_k(x)$
the upper tail probability of $\dgamma(k\shape,1)$, i.e.
$G_k(x) = \nth{\Gamma(k\shape)}\int_x^\infty s^{k\shape-1}
e^{-s}\,ds$, $x>0$.  As noted in Section \ref{subsec:gamma}, the main
step is to find the asymptotics of $G_k(c_k)$ and $G_k(d_k)$, where
$c_k$ is defined in \eqref{eq:gamma-thresh} and
$d_k = \frac{c_k}{1+\delta}$.

To find the asymptotic of $G_k(c_k)$, first, by power expansion of
$\ln(1+\delta)$, for $\delta\in (-1,1)$, 
\begin{align*}
  c_k
  &=
  k\shape + \frac{\ln Q_\alpha}{\delta} + k\shape\delta
  \sum_{j=0}^\infty \frac{(-\delta)^j}{(j+1)(j+2)} + \ln
  Q_\alpha, \\
  d_k 
  &= k\shape + \frac{\ln Q_\alpha}{\delta}
  + k\shape\delta 
  \sum_{j=0}^\infty \frac{(-1)^{j-1}\delta^j}{j+2}.
\end{align*}
Because $k\delta\toi$ while $k\delta^2\to 0$, it is seen that in each
of the sums, every term is of an infinitesimal order of its previous
one.  Let
\begin{align}
  z = k\shape,\ \ s = z(1+t), \ \ b = \ln Q_\alpha.
  \label{eq:zsb}
\end{align}
With the variable substitutions,
\begin{align*}
  G_k(c_k)
  = \frac{z^z e^{-z}}{\Gamma(z)} 
  \underbrace{\int_{\frac{b}{z\delta} + D(\delta, z)}^\infty
    (1+t)^{z-1} e^{-z t} \,d t}_{I(z,\delta)}, \quad
  G_k(d_k)
  = \frac{z^z e^{-z}}{\Gamma(z)} 
  \underbrace{
    \int_{\frac{b}{z\delta} + L(\delta, z)}^\infty(1+t)^{z-1} e^{-z t}
    \,d t}_{J(z,\delta)},
\end{align*}
where
\begin{align}
  D(\delta, z)& = \delta r(\delta) + \frac{b}{z},
  \text{ with }
  r(\delta) = \sum_{j=0}^\infty
  \frac{(-\delta)^j}{(j+1)(j+2)}, \label{eq:D} \\
  L(\delta, z)
  &= \delta \bar r(\delta) + \frac{b}{z}, 
  \text{ with }
  \bar r(\delta) = -\sum_{j=0}^\infty \frac{(-\delta)^j}{j+2}.
  \label{eq:L}
\end{align}
The main step is to show
\begin{align}
  I(z,\delta)
  &
  \sim \frac{\delta}{b}
  \exp\Cbr{
    -\frac{b^2}{2z\delta^2} -\frac{b}{2} - z\psi\Grp{
      \frac{b}{z\delta}}
  }, \label{eq:gamma-I} \\
  J(z,\delta)
  &
  \sim
  \frac{\delta}{b}
  \exp\Cbr{
    -\frac{b^2}{2z\delta^2} + \frac{b}{2} - z\psi\Grp{
      \frac{b}{z\delta}}
  }.
  \label{eq:gamma-J}
\end{align}

Assume the two formulas are true for now.  By Stirling's formula, 
$\frac{z^z e^{-z}}{\Gamma(z)}
= (1+o(1))\sqrt{\frac{z}{2\pi}}$.  Then by \eqref{eq:gamma-I} and
\eqref{eq:gamma-J},
\begin{align*}
  G_k(c_k)
  &
  \sim \sqrt{\frac{z}{2\pi}}
  \frac{\delta}{b}
  \exp\Cbr{
    -\frac{b^2}{2z\delta^2} -\frac{b}{2} - z\psi\Grp{
      \frac{b}{z\delta}}
  }, \\
  G_k(d_k)
  &\sim \sqrt{\frac{z}{2\pi}}
  \frac{\delta}{b}
  \exp\Cbr{
    -\frac{b^2}{2z\delta^2} + \frac{b}{2} - z\psi\Grp{
      \frac{b}{z\delta}}
  }.
\end{align*}
Since $p_{k,\delta}(\alpha) = (1-\pa) G_k(c_k)+\pa G_k(d_k)$ and
\eqref{eq:zsb},
\begin{align*}
  p_{k,\delta}(\alpha)
  \sim\
  &
  \sqrt{\frac{(1-\pa)\pa}{2\pi(1-\alpha)\alpha}}
  \frac{\sqrt{k\shape} \delta}{\sqrt{\ln Q_\alpha}}
  \exp\Cbr{
    -\frac{(\ln Q_\alpha) ^2}{2k\shape \delta^2}
    - k\shape \psi\Grp{
      \frac{\ln Q_\alpha}{k\shape\delta}}
  }.
\end{align*}
The proof is complete by \eqref{eq:min-n}.

The rest of the section is devoted to the proof of \eqref{eq:gamma-I}
and \eqref{eq:gamma-J}.  Observe $r(\delta)\to 1/2$ and $\bar
r(\delta)\to -1/2$ as $\delta\to 0$.  By \eqref{eq:k-delta},
\begin{align}
  z\delta = \shape k\delta \toi,  \ \ z\delta^2 = \shape k\delta^2 \to
  0.
  \label{eq:z-delta}
\end{align}
It is then not hard to check that $D(\delta, z)\sim \delta/2$.  Also,
for any $z>0$, $(1+t)^{z-1} e^{-z t}$ is strictly decreasing in $t>0$.
Given $0<\epsilon \ll 1$, using \eqref{eq:z-delta} again, 
\begin{align*}
  I_\epsilon(z,\delta)&:=
  \int_{\frac{b}{z\delta}+D(\delta,z)}^\epsilon
  (1+t)^{z-1} e^{-z t}\,d t
  \ge
  \int_{\frac{2b}{z\delta}}^{\frac{3b}{z\delta}}
  (1+t)^{z-1} e^{-z t}\,d t \\
  &
  \ge \Grp{1+\frac{3b}{z\delta}}^{z-1} e^{-3b/\delta}\Grp{
    \frac{b}{z\delta}
  } \\
  &
  \stackrel{(a)}{\ge}
  \exp\Cbr{
    \Sbr{
      \frac{3b}{z\delta} - \nth 2
      \Grp{\frac{3b}{z\delta}}^2
    }(z-1) - \frac{3b}{\delta}
  } \frac{b}{z\delta} \ge
  \exp\Cbr{- \frac{C}{z\delta^2}} \frac{b}{z\delta}
\end{align*}
for some $C>0$, where $(a)$ is due to $\ln(1+x) \ge x - x^2/2$ for
$x>0$.  Therefore, $I_\epsilon(z, \delta)^{1/z}\to 1$.  On the other hand,
$$
\Grp{
  \int_\epsilon^\infty (1+t)^{z-1} e^{-z t}\,d t
}^{1/z} \to \sup_{t\ge\epsilon} (1+t) e^{-t} < 1.
$$
As a result, for any $\epsilon>0$, $I(z,\delta) \sim I_\epsilon(z,
\delta)$.  Since $\epsilon$ is arbitrary, it follows that we can
replace $(1+t)^{z-1}$ in the integrand to $(1+t)^z$ to get
\begin{align}
  I(z,\delta) \sim \int_{\frac{b}{2\delta} + D(\delta, z)}^\epsilon
  (1+t)^z e^{-z t}\,d t =
  \int_{\frac{b}{2\delta} + D(\delta, z)}^\epsilon
  e^{-z\varphi(t)}\,d t,  \label{eq:I-phi}
\end{align}
where $\varphi(t)= t - \ln(1+t)$.  By $\varphi'(t) = \frac{t}{1+t}$,
$\varphi(t)$ is a strictly increasing function from $(0,\infty)$ onto
$(0,\infty)$ with smooth inverse $\varphi^{-1}(u)$.  On the other
hand, $\varphi(t) = \frac{t^2}{2} + O(t^3)$, as $t\downarrow 0$.  As a
result, as $u\to 0+$, $\varphi^{-1}(u) = (1+o(1)) \sqrt{2u}$ and hence
$$
(\varphi^{-1})'(u) = \nth{\varphi'(\varphi^{-1}(u))}
= 1+\nth{\varphi^{-1}(u)} = \frac{1+o(1)}{\sqrt{2u}} \quad
\text{ as } u\to 0.
$$
By \eqref{eq:I-phi} and the arbitrariness of $\epsilon>0$ as well as
the above properties of $\varphi$,
\begin{align*}
  I(z,\delta)
  \sim
  \int_{\varphi(\frac{b}{z\delta}+D(\delta,z))}^{\varphi(\epsilon)}
  \frac{e^{-z u}}{\sqrt{2u}}\,d u 
  \sim
  \int_{\varphi(\frac{b}{z\delta}+D(\delta,z))}^\infty 
  \frac{e^{-z u}}{\sqrt{2u}}\,d u  = I_1.
\end{align*}
By variable substitution $u = v/z$,
\begin{align*}
  I_1 = \nth{\sqrt{z}}
  \int_{z\varphi(\frac{b}{z\delta}+D(\delta,z))}^\infty
  \nth{\sqrt{2v}} e^{-v}\,d v.
\end{align*}
Since $z\delta D(\delta,z)\toi$ and $D(\delta,z)\to 0$, by
$\varphi(t)\sim t^2/2$ as $t\to 0$, 
\begin{align}
  z\varphi\Grp{\frac{b}{z\delta} + D(\delta,z)} \sim
  \frac{z}{2} \Grp{\frac{b}{z\delta}}^2 = \frac{b^2}{2 z\delta^2}
  \toi.
  \label{eq:z-phi}
\end{align}
Recall that for any $a$, $\int_x^\infty t^a e^{-t}\,d t \sim x^a
e^{-x}$, as $x\toi$.  Then by \eqref{eq:z-phi}
\begin{align*}
  I_1
  &\sim \nth{\sqrt{z}}
  \nth{\displaystyle
    \sqrt{2 z \varphi\Grp{\frac{b}{z\delta} + D(\delta,z)}
    }
  } \exp\Cbr{-z \varphi\Grp{\frac{b}{z\delta}+D(\delta,z)}} \\
  &
  \sim \frac{\delta}{b} \exp\Cbr{
    -z \varphi\Grp{\frac{b}{z\delta}+D(\delta,z)}
  }.
\end{align*}
Because $\psi(t) = \varphi(t) - \frac{t^2}{2}$,
\begin{align*}
  z\varphi\Grp{\frac{b}{z\delta}+D(\delta,z)}
  = \frac{z}{2}\Grp{\frac{b}{z\delta}+D(\delta,z)}^2 +
  z\psi\Grp{\frac{b}{z\delta}+D(\delta,z)}.
\end{align*}
First, by \eqref{eq:z-delta},
\begin{align*}
  &\hspace{-1cm}
  \frac{z}{2}\Grp{\frac{b}{z\delta}+D(\delta,z)}^2 
  = \frac{b^2}{2z\delta^2} + \frac{b D(\delta,z)}{\delta}
  + \frac{z(D(\delta,z))^2}{2} \\
  &
  =
  \frac{b^2}{2z\delta^2} + \frac{b}{\delta}
  \Grp{\delta r(\delta)+\frac{b}{z}}
  + \frac{z}{2}\Grp{\delta r(\delta)+\frac{b}{z}}^2 
  = \frac{b^2}{2z\delta^2} + \frac{b}{2} + o(1).
\end{align*}
Second, since $\psi'(x) = \varphi'(x) - x = 1-\frac{1}{1+x}-x = -
\frac{x^2}{1+x}$, by Taylor expansion and \eqref{eq:D}
\begin{align*}
  z\psi\Grp{\frac{b}{z\delta}+D(\delta,z)}-
  z\psi\Grp{\frac{b}{z\delta}}
  =
  z D(\delta,z) \psi'\Grp{\frac{b}{z\delta}+\xi D(\delta,z)}
  =
  \frac{z\delta}{2}\Grp{\frac{b}{z\delta}}^2 = o(1).
\end{align*}
As a result,
\begin{align*}
  z\varphi\Grp{\frac{b}{z\delta}+D(\delta,z)}
  = \frac{b^2}{2z\delta^2} + \frac{b}{2} + z\psi\Grp{
    \frac{b}{z\delta}} + o(1)
\end{align*}
and hence by $I(z,\delta)\sim I_1$, \eqref{eq:gamma-I} then follows.
By similar argument, it can be shown that
\begin{align*}
  J(z,\delta)
  &\sim
  \nth{\sqrt{z}}
  \nth{\displaystyle
    \sqrt{2 z \varphi\Grp{\frac{b}{z\delta} + L(\delta,z)}
    }
  } \exp\Cbr{-z \varphi\Grp{\frac{b}{z\delta}+L(\delta,z)}} \\
  &
  \sim \frac{\delta}{b} \exp\Cbr{
    -z \varphi\Grp{\frac{b}{z\delta}+L(\delta,z)}
  },
\end{align*}
which leads to \eqref{eq:gamma-J}.

\subsection{Proofs for the asymptotics of power}
A basic fact to use is that under the random effects
model \eqref{eq:model}, $\Pr(\eta_i=0\gv\vf X)$ are i.i.d.\ and by
\eqref{eq:posterior}, given $\delta>0$ and $k$, for any $\alpha\in
(0,1)$, the probability of $\{\Pr(\eta_i=0\gv \vf X)\le\alpha\}$ is
\begin{align}
  p_{k,\delta}(\alpha)=(1-\pa)
  P_0(E_k(\alpha)) + \pa P_a(E_k(\alpha)),
  \label{eq:eta}
\end{align}
where $P_0$ and $P_a$ are the probability distributions under true and
false nulls, respectively, and $E_k(\alpha) = \Cbr{\prod_{j=1}^k
  \frac{f_a(X_j)}{f_0(X_j)} \ge Q_\alpha }$ with $\eno X k$ being
i.i.d. 
\paragraph{\it Proof of Proposition \ref{prop:odds-ratio}.}
For fixed $N$,
\begin{align*}
  R = \sum_{i=1}^N d^*_i(\vf X), \quad
  R_0 = \sum_{i=1}^N d^*_i(\vf X)(1-\eta_i), \quad
  R_a = \sum_{i=1}^N d^*_i(\vf X)\eta_i.
\end{align*}
Given $\delta>0$ and $k$, since $d_i^*(\vf X) = \cf{\Pr(\eta_i=0\gv
  \vf X)\le \alpha}$, by \eqref{eq:eta} and the Weak Law of Large
Numbers (WLLN), $R/N\convP p_{k,\delta}(\alpha)>0$.  Similarly,
$R_0/N\convP (1-a)P_0(E_k(\alpha))$, $R_a/N\convP a P_a(E_k(\alpha))$
and $N_a/N\to a$.  Property \eqref{eq:property} can then be proved.
In particular,
\begin{align}
  \power_\infty(\delta, k) = P_a(E_k(\alpha)), \quad
  \pfdr_\infty^*(\delta,k) = \frac{(1-a)
    P_0(E_k(\alpha))}{p_{k,\delta}(\alpha)}.
  \label{eq:finite}
\end{align}

To show \eqref{eq:odds-large-n}, given $\vf X$ with $R>0$,
\begin{align*}
  E\Sbr{R_0/R\gv \vf X}
  &
  = \nth{R} \sum_{i=1}^N 
  E\Sbr{
    \cf{\Pr(\eta_i=0\gv\vf X)\le\alpha} (1-\eta_i)\gv \vf X
  } \\
  &
  = \nth{R} \sum_{i=1}^N 
  \cf{\Pr(\eta_i=0\gv\vf X)\le\alpha} \Pr(\eta_i=0\gv\vf X) \\
  &
  \le \nth{R} \sum_{i=1}^N 
  \alpha \cf{\Pr(\eta_i=0\gv\vf X)\le\alpha} = \alpha.
\end{align*}
Since $p_{k,\delta}(\alpha)>0$, $\Pr(d^*_i(\vf X)>0$ for at least one
$i=1, \ldots, N)>0$.  Therefore, the conditional expectation of $R_0/R$
over $\vf X$ with $R>0$ is well defined, giving $E[R_0/R\gv R>0]
\le\alpha$.  Thus $\pfdr^*_\infty(\delta,k)\le \alpha$ for all
$(\delta, k)$.  On the other hand, given $\beta<\alpha$,
\begin{align*}
  E[R_0/R \gv \vf X] 
  &
  \ge
  \nth R\sum_{i=1}^N E[\cf{\beta\le \Pr(\eta_i=0\gv \vf X) \le \alpha}
  (1-\eta_i) \gv \vf X] \\
  &
  \ge 
  \frac{\beta}{R}\sum_{i=1}^N \cf{\beta\le \Pr(\eta_i=0\gv \vf X) \le
    \alpha}.
\end{align*}
Under the random effects model and the WLLN,
\begin{gather*}
  \sum_{i=1}^N \cf{\beta\le \Pr(\eta_i=0\gv \vf X)\le \alpha}
  = (1+o_p(1)) [p_{k,\delta}(\alpha) - p_{k,\delta}(\beta)]N,
\end{gather*}
where $o_p(1)$ stands for some sequence of random variables $\xi_N
\convP 0$ as $N\toi$.  Taking expectation over $\vf X$ with $R>0$ 
and then letting $N\toi$,
$$
\pfdr^*_\infty(\delta,k) \ge\frac{\beta[p_{k,\delta}(\alpha) -
  p_{k,\delta}(\beta)]}{p_{k,\delta}(\alpha)}.
$$
Let $(\delta, k)\to (0,\infty)$ while satisfying \eqref{eq:positive}
and \eqref{eq:increase}.  Then $\Limsup \pfdr^*_\infty(\delta, k) \ge
\beta$.  Since $\beta$ is arbitrary, $\pfdr^*_\infty(\delta,k) \to
\alpha$, showing the second half of \eqref{eq:odds-large-n}.
Finally, combining this with \eqref{eq:eta} and \eqref{eq:finite}, 
the first half of \eqref{eq:odds-large-n} follows.  \qed

\paragraph{\it Proof of Proposition \ref{prop:large-n}.} 
Let $d$ be a procedure more powerful than $d^*$ while satisfying
\eqref{eq:property} and $\Limsup\pfdr_\infty(\delta,k)\le \alpha$.
Let $(\delta, k)$ be fixed first.  Given $0<\alpha_1 <
\alpha<\alpha_2$, let
\begin{align*}
  R\Sp 1
  &
  = \#\{i: \Pr(\eta_i=0\gv\vf X)<\alpha_1,\ d_i(\vf X)=1\},
  \\
  R\Sp 2
  &
  = \#\{i: \alpha_1\le \Pr(\eta_i=0\gv\vf X)\le \alpha_2,\
  d_i(\vf X)=1\},
  \\
  R\Sp 3
  &
  = \#\{i: \Pr(\eta_i=0\gv\vf X)> \alpha_2,\
  d_i(\vf X)=1\}.
\end{align*}
Then $R=R\Sp 1+R\Sp 2+R\Sp 3$ and for any
$\vf X$ with $R>0$,
\begin{align}
  E\Sbr{R_0/R \gv \vf X} = \frac{1}{R}
  \sum_{i=1}^N d_i(\vf X) \Pr(\eta_i=0\gv \vf X)
  \ge \frac{\alpha_1 R\Sp 2 + \alpha_2 R\Sp 3}{R}.
  \label{eq:pfdr-X}
\end{align}
Since $d$ satisfies \eqref{eq:property}, by the WLLN, $R_a =
(\pa+o_p(1))\power_\infty(\delta,k) N$ as $N\toi$.  By the assumption
and \eqref{eq:odds-large-n}, $\power_\infty(\delta,k)\ge
\power_\infty^*(\delta,k)=(1-\alpha + o(1)) p_{k,\delta}(\alpha)/\pa$.
Since $R\ge R_a$, $R$ is at least of the same order as
$p_{k,\delta}(\alpha)N$.  On the other hand, since
$$
R\Sp 1 \le \#\{i: \Pr(\eta_i=0\gv\vf X)<\alpha_1\} = (1+o_p(1))
p_{k,\delta}(\alpha_1) N,
$$
by \eqref{eq:increase}, $R\Sp 1 = o_p(1)R$ and so $R\Sp 2+R\Sp 3=
(1+o_p(1)) R$.  Therefore,
\begin{align*}
  \frac{\alpha_1 R\Sp 2 + \alpha_2 R\Sp 3}{R}
  =(1+o_p(1)) 
  \frac{\alpha_1 R\Sp 2 + \alpha_2 R\Sp 3}{R\Sp 2+R\Sp 3}
  = (1+o_p(1)) \Sbr{\alpha_1 + (\alpha_2-\alpha_1) \frac{R\Sp 3}{R}}.
\end{align*}
Combine this with \eqref{eq:pfdr-X}.  Taking expectation over $\vf X$
and letting $N\toi$ yield
\begin{align*}
  \pfdr_\infty(\delta,k)
  \ge \alpha_1 + (\alpha_2 - \alpha_1) \Limsup_{N\toi}
  E[\,R\Sp 3/R\gv R>0\,].
\end{align*}

Let $(\delta,k)\to (0,\infty)$.  As $\Limsup\pfdr_\infty(\delta, k)\le
\alpha$ by the assumption and $\alpha_1<\alpha$ is arbitrary, the
second half of \eqref{eq:large-n} follows.  Furthermore, the above
inequality implies $\Limsup_N E[R\Sp 3/R\gv R>0] = o(1)$.  Since for
each fixed $(\delta, k)$, $R\Sp 1 = o_p(1) R$ and $N_a = (\pa +o_p(1)) N$,
it then follows that, as $(\delta,k)\to (0,\infty)$, $\lim_N E[R/N_a]
= \lim_N E[(R\Sp 2 + R\Sp 3)/N_a] = (1+o(1)) \lim_N E[R\Sp
2/N_a]$ .  Since $R\Sp 2$ is no greater
than the number of nulls with $\Pr(\eta_i=0\gv\vf X)\le \alpha_2$,
which is $(1+o_p(1))p_{k,\delta}(\alpha_2) N$,
\begin{align*}
  \power_\infty(\delta, k)\le\lim_{N\toi} E\Sbr{\frac{R}{N_a}}
  =
  (1+o(1)) \lim_{N\toi} E\Sbr{\frac{R\Sp 2}{N_a}}
  \le (1+o(1)) \frac{p_{k,\delta}(\alpha_2)}{\pa}.
\end{align*}
Therefore, $\power_\infty(\delta,k)$ satisfies \eqref{eq:large-n}.
\qed

\paragraph{Proof of Proposition \ref{prop:power}.}  Denote $\alpha_k =
\alpha + c k\delta^2$.  By the WLLN, $d$ satisfies
\eqref{eq:property}.  For each $(\delta,k)$, $d$ is a 
thresholding procedure, so $\pfdr_\infty(\delta,k) \le
\alpha_k$.  Then by $k\delta^2 \to 0$, $d$ satisfies $\Limsup
\pfdr_\infty(\delta,k)\le \alpha$.  Given $(\delta,k)$,
following the same argument that leads to \eqref{eq:p-lrt1}, 
\begin{align*}
  \power_\infty(\delta,k) = 
  (1+o(1)) 
  \frac{\sqrt{ k I(\theta_0)} \delta}{\sqrt{2\pi}
    \ln Q_{\alpha_k}}
  \exp\Cbr{
    -\frac{(\ln Q_{\alpha_k})^2}{2 k I(\theta_0) \delta^2}
    + \frac{\ln Q_{\alpha_k}}{2}
  }.
\end{align*}
Since $Q_{\alpha_k}\to Q_\alpha$, to get $\power_\infty(\delta,k) =
(M+o(1)) \power_\infty^*(\delta,k)$ as $(\delta,k)\to (0,\infty)$, it
boils down to showing
$$
\exp\Cbr{
  -\frac{(\ln Q_{\alpha_k})^2}{2 k I(\theta_0) \delta^2}
} = (1+o(1)) M
\exp\Cbr{
  -\frac{(\ln Q_\alpha)^2}{2 k I(\theta_0) \delta^2}
}.
$$
By Taylor expansion,
$$
\ln Q_{\alpha_k} - \ln Q_\alpha 
= \ln \frac{1-\alpha_k}{1-\alpha} - \ln \frac{\alpha_k}{\alpha}
= - \frac{c k\delta^2}{(1-\alpha)\alpha} + O((k\delta^2)^2).
$$
As a result,
$$
\frac{(\ln Q_{\alpha_k})^2}{2kI(\theta_0)\delta^2}
= 
\frac{(\ln Q_\alpha)^2}{2kI(\theta_0)\delta^2}
- \frac{c\ln Q_\alpha}{(1-\alpha)\alpha I(\theta_0)} + O(k\delta^2).
$$
By the definition of $c$, the result follows. 

Finally, we show that among procedures that satisfy
\eqref{eq:property} and $\Limsup\pfdr_\infty(\delta,k)\le \alpha$, no
one is asymptotically the most powerful.   It suffices to show that
for any such procedure $d$ that is more powerful than $d^*$, there is
yet another one more powerful than $d$.  First, by diagonal argument
and the first part of \eqref{eq:large-n}, there is a decreasing
$\alpha_k\to\alpha$, such that $\power_\infty(\delta, k) \le 2
p_{k,\delta}(\alpha_k)/\alpha$ for large $k$.  Now we use the same
construction as above.  Let $\alpha_k' = \alpha_k + c_k k\delta^2$,
with $c_k = (1-\alpha_k) \alpha_k I(\theta_0) M/\ln Q_{\alpha_k}$,
where $M>0$.  The thresholding procedure using $\alpha_k'$ as cut-offs
satisfy the conditions of Proposition \ref{prop:large-n}.  It is seen
that as long as $M$ is large enough, the power of this new procedure
will be greater than $4 p_{k,\delta}(\alpha_k)/\alpha$, and hence at
least twice as large as $\power_\infty(\delta,k)$.
\qed

\end{document}